\documentclass{amsart}

\usepackage{graphicx}
\usepackage{latexsym}
\usepackage{amsfonts}
\usepackage[all]{xy}
\usepackage{amssymb, mathrsfs, amsfonts, amsmath}
\usepackage{amsbsy}
\usepackage{fontenc}
\usepackage{amssymb}
\newtheorem{theorem}{Theorem}[section]

\newtheorem{conjecture}[theorem]{Conjecture}

\theoremstyle{definition}

\theoremstyle{remark}

\numberwithin{equation}{section}

\newcommand{\diver}{{{\rm{div}\,}}}

\newcommand{\hess}{\textrm{Hess}}
\newcommand{\grad}{\textrm{grad}\,}

\newcommand{\noin}{\noindent}

\begin{document}

\title{On the essential   spectrum of  Nadirashvili-Martin-Morales minimal surfaces}

\author{G. Pacelli Bessa}
\curraddr{Department of Mathematics, Universidade Federal do Ceara-UFC, Campus do Pici, 60455-760 Fortaleza-CE
Brazil}
\email{bessa@mat.ufc.br}
\thanks{The first author was partially supported  by  a CNPq-grant}

\author{Luquesio P. Jorge}
\email{ljorge@mat.ufc.br}
\author{J. Fabio Montenegro}
\email{fabio@mat.ufc.br}

\date{\today }

\dedicatory{}


\begin{abstract}
  We show  that the spectrum of a complete submanifold properly immersed into a ball of a Riemannian manifold is discrete, provided the norm of the mean curvature vector is sufficiently small. In particular,  the spectrum of a  complete minimal surface properly immersed into a ball of $\mathbb{R}^{3}$ is  discrete. This gives a positive answer to a question of Yau \cite{Yau3}.
\end{abstract}

\maketitle

\section{Introduction}\label{section1}
\label{intro}An interesting problem in the Geometry of the Laplacian  is to understand the relations of the geometry of a Riemannian manifold and its spectrum. For instance, to understand the
restrictions on the geometry of a Riemannian manifold implying  that its spectrum  is purely continuous or discrete. There are several important work along these lines. See \cite{donnelly1}, \cite{donnelly-garofalo}, \cite{escobar}, \cite{karp}, \cite{rellich}, \cite{tayoshi} for geometric conditions implying that the spectrum is purely continuous and  \cite{baider},  \cite{donnelly-li}, \cite{harmer}, \cite{kleine1}, \cite{kleine2}  for geometric conditions implying  that the spectrum is discrete.

Since
every  complete Riemannian $m$-manifold can be realized as a complete submanifold embedded into a ball of radius $r$ of an $n$-dimensional Euclidean  space, with $n$ depending only on $m$, see \cite{nash}, it would be  important to understand  the relations between  the spectrum and the extrinsic geometry of bounded embeddings of complete Riemannian manifolds in Euclidean spaces.   A particularly interesting aspect of this problem is the spectrum related  part of the so called {\em Calabi-Yau conjectures on minimal surfaces}.

Yau in his $2000$ millennium lectures  \cite{Yau3}, \cite{yau4}, revisiting these conjectures, wrote:
{\em  It is known \cite{Nadirashvili} that there are complete minimal surfaces properly immersed into a $[$open$]$ ball. $\ldots$ Are their spectrum discrete?} It is worthwhile to point out that it is not clear that the  Nadirashvili's complete bounded minimal surface \cite{Nadirashvili} is {\em properly immersed}. However, in \cite{pacomartin},  \cite{pacomartin2}, F.  Martin and S. Morales   constructed, for any open convex subset $B$ of $ \mathbb{R}^{3}$, a complete proper minimal immersions $\varphi \colon \mathbb{D} \hookrightarrow B$, where $\mathbb{D}$ is the standard disk on $\mathbb{R}^{2}$. The Martin-Morales' method is a highly non-trivial refinement of  Nadirashvili's method, thus we  name, (as we should),  these complete properly immersed minimal surfaces into  convex subsets $B$ of $\mathbb{R}^{3}$  as  Nadirashvili-Martin-Morales minimal surfaces.

The purpose of this paper is to  answer positively Yau's question. In fact,  we show as a particular case of our main result that the spectrum of  any  Nadirashvili-Martin-Morales minimal surface  is discrete if the convex set $B$ is a ball $B_{\mathbb{R}^{3}}(r)$ of $\mathbb{R}^{3}$. We prove the following.
  \begin{theorem} \label{thm1}Let $\varphi : M \hookrightarrow B_{\mathbb{R}^{3}}(r)\subset \mathbb{R}^{3}$ be a  complete   surface, properly immersed into a ball.  If the norm of the  mean curvature vector $H$ of $M$ satisfies $$\sup_{M} \vert H\vert < 2/r$$ then  $M$ has discrete spectrum.
\end{theorem}

Our main result Theorem \ref{thmMain}\hspace{1mm} is a natural generalization of Theorem \ref{thm1}. It shows that the spectrum of a complete properly immersed submanifold $\varphi \colon M\hookrightarrow B_{N}(r)$ is discrete provided the norm of the mean curvature vector  $ H={\rm Tr}\,\alpha $ is sufficiently small. Here $B_{N}(r)\subset N$ is a normal geodesic ball of radius $r$ of  a Riemannian manifold $N$ and $\alpha $ is the second fundamental form.
In the following we denote
\begin{equation}\label{eqCb}
C_b(t)
=\left\{\begin{array}{lll}
\sqrt{b}\cot(\sqrt{b}\, t) & \mathrm{if} & b >0,\\
1/t & \mathrm{if} & b =0,\\
\sqrt{-b}\coth(\sqrt{-b}\, t) & \mathrm{if}  & b <0.
\end{array}\right.
\end{equation}

\begin{theorem}\label{thmMain}
Let $\varphi\colon M\hookrightarrow  B_{N}(r)$ be a  complete   $m$-submanifold  properly immersed into a geodesic ball, centered at $p$ with radius $r$, of a  Riemannian $n$-manifold $N$.  Let  $b=\sup K_N^{\mathrm{rad}}$ where $K_N^{\mathrm{rad}} $ are  the radial sectional
curvatures along the  geodesics issuing from
$p$. Assume that  $r<\min\{{\rm inj}_N(p), \pi/2\sqrt{b}\}$, where $\pi/2\sqrt{b}=+\infty$ if $b\leq 0$.
If If the norm of the  mean curvature vector $H$ satisfies, $$
\sup_{M}\vert  H\vert <m\cdot C_{b}(r),
$$
then $M$ has discrete spectrum.
\end{theorem}

The properness condition  is a marginal technical hypothesis  in Theorem (\ref{thmMain}).
 It is used only to choose a natural sequence of compact subsets of $M$ so that we can construct a sequence of positive smooth functions on their complements. The result  should hold without it.

Isabel Salavessa  in a beautiful paper \cite{salavessa}, generalized Theorem (\ref{thmMain}) in the minimal case proving discreteness of the spectrum of  $X$-bounded minimal submanifolds of  Riemannian manifolds carrying strongly convex vector field $X$.

A Riemannian
manifold $M$ is said to be \emph{stochastically complete} if for some (and therefore, for any) $(x,t)\in M\times(0,+\infty)$ it holds
that
$$
\int_Mp(x,y,t)dy=1,
$$
where $p(x,y,t)$ is the heat kernel of the Laplacian operator. Otherwise, the manifold $M$ is said to
be \emph{stochastically incomplete} (for further details about this see, for instance, \cite{Gr}).
It seems to have a close relation between discreteness of the spectrum of a complete noncompact Riemannian manifolds and stochastic incompleteness. For instance, it  was proved in \cite{alias-bessa-dajczer} that  submanifolds satisfying the hypotheses of Theorem (\ref{thmMain}), (without the properness condition)  are stochastically incomplete. M. Harmer \cite{harmer}, shows that stochastic incompleteness implies discreteness of the spectrum in a certain class of Riemannian manifolds. Based on these evidences, we believe that  the following conjecture should be true.
\begin{conjecture}A complete noncompact Riemannian manifold has discrete spectrum if and only if is stochastically incomplete.
\end{conjecture}

\section{ Preliminaries.}\label{section2}Let $M$ be a  complete noncompact Riemannian manifold. The Laplacian $\triangle$ acting on $C_{0}^{\infty}(M)$ has a  unique self-adjoint extension
 to an unbounded operator acting
on $L^{2}(M)$, also  denoted by $\triangle$, whose domain are those functions  $f\in L^{2}(M)$ such that $\triangle f\in L^{2}(M)$ and whose spectrum  $\Sigma(M)\subset [0,\infty)$
decomposes as $\Sigma(M)=\Sigma_{p}(M)\cup\Sigma_{ess}(M)$ where $\Sigma_{p}(M)$
is  formed by   eigenvalues  with finite multiplicity
 and $\Sigma_{ess}(M)$ is  formed
by accumulation points of the spectrum and by the eigenvalues  with
infinite multiplicity. It is said that
 $M$ has discrete spectrum if $\Sigma_{ess}(M)=\emptyset$
  and that  $M$ has purely continuous spectrum if  $\Sigma_{p}(M)=\emptyset$.

 If $K\subset M$ is a compact manifold with boundary, of the same dimension as $M$ then there is a self-adjoint extension $\triangle'$ of the Laplacian $\triangle$ of $M\setminus K$ by imposing Dirichlet conditions.
The {\em Decomposition Principle} \cite{donnelly-li} says that $\triangle$ and $\triangle'$ have the same essential spectrum $\Sigma_{ess}(M)=\Sigma_{ess}(M\setminus K)$. On the other hand, the bottom of the spectrum of $M\setminus K $   is equal to the fundamental tone of $M\setminus K$, i.e.    $\inf \Sigma (M\setminus K)=\lambda^{\ast}(M\setminus K)$, where
$$\lambda^{\ast}(M\setminus K)=\inf\left\{\frac{\int_{M\setminus K}\vert \grad f\vert^{2}}{\int_{M\setminus K}f^{2}}, f\in C_{0}^{\infty}(M\setminus K)\setminus \{0\}\right\}.$$To give lower estimates for $\lambda^{\ast}(M\setminus K)$ we need of the following version of Barta's Theorem.
\begin{theorem}[Barta, \cite{barta}]\label{barta}
 Let $\Omega\subset M$ be an open subset of  a Riemannian manifold $M$ and let
$f\in C^{2}(\Omega)$, $f\vert\Omega >0$. Then
\begin{equation}\label{eq1.3}
\lambda^{\ast}(\Omega ) \geq \inf_{\Omega} (-\Delta f/f).
\end{equation}
\end{theorem}
\noin {\em Proof:}  Let $X=-\grad \log f$ be a $C^{1}$ vector field in $\Omega$. It was proved in \cite{bessa-montenegro2} that $$\lambda^{\ast}(\Omega)\geq \displaystyle \inf_{\Omega} (\diver X - \vert X\vert^{2})=\inf_{\Omega}( -\frac{\triangle f}{f}).$$The second main ingredient of our proof is the Hessian comparison theorem.

\begin{theorem}\label{thm2}
Let $M^m$ be a  Riemannian manifold  and
$x_0,x_1 \in M$ be such that there is a minimizing unit speed geodesic $\gamma$  joining $x_{0}$ and $x_{1}$ and let $\rho(x)=\mathrm{dist}(x_0,x)$
be the  distance function  to $x_{0}$. Let $a\leq K_{\gamma}\leq b$ be the radial
sectional curvatures of $M$ along $\gamma$.
If $b>0$ assume $\rho(x_{1})<\pi/2\sqrt{b}$. Then, we have $\hess\rho(x) (\gamma',\gamma')=0$ and
\begin{equation}\label{thmHess}
C_{a}(\rho(x))\Vert
X\Vert^2\geq \hess\rho(x)(X,X)\geq C_{b}(\rho(x))\Vert
X\Vert^2
\end{equation}
where $X\in T_{x}M$ is perpendicular to $\gamma'(\rho(x))$.
\end{theorem}

 Let $\varphi : M \hookrightarrow W$ be an isometric  immersion of a complete  Riemannian $m$-manifold $M$ into a Riemannian $n$-manifold $W$  with second fundamental form $\alpha$.  Consider a $C^{2}$-function $g:W \rightarrow \mathbb{R}$ and the composition $f=g\,\circ\,
 \varphi :M \rightarrow \mathbb{R}$. Identifying $X$ with $d\varphi (X)$ we have at $q\in M$ that
 the Hessian of $f$ is given by
 \begin{equation}\hess\,f (q)  \,(X,Y)= \hess\,g (\varphi (q))\,(X,Y) +
 \langle \grad\,g\,,\,\alpha (X,Y)\rangle_{\varphi (q)}.
\label{eqBF2}
\end{equation}
 Taking the trace in (\ref{eqBF2}), with respect to an orthonormal basis $\{ e_{1},\ldots e_{m}\}$
 for $T_{q}M$, we have the Laplacian of $f$,
\begin{equation}\begin{array}{lll}
\Delta \,f (q) & = & \sum_{i=1}^{m}\hess\,g (\varphi
(q))\,(e_{i},e_{i}) + \langle \grad\,g\,,\, \sum_{i=1}^{m}\alpha
(e_{i},e_{i})\rangle.\label{eqBF3}\end{array}
\end{equation}
The formulas (\ref{eqBF2}) and (\ref{eqBF3}) are  well known in
the literature,
 see  \cite{jorge-koutrofiotis}.

 \section{Proof of Theorem \ref{thmMain}}\label{subsection3} Let $K_{1}\subset K_{2}\subset \cdots $ be an exhaustion sequence of $M$ by compact sets. The {\em Decomposition Principle} states that $M$ and $M\setminus K_{i}$ have the same essential spectrum, $\Sigma_{ess}(M)=\Sigma_{ess}(M\setminus K_{i})$. Therefore, the Theorem \ref{thmMain} is proved if we show that $\lim_{i\to \infty}\lambda^{\ast}(M\setminus K_{i})=\infty$ since $\lambda^{\ast}(M\setminus K_{i})\leq \inf \Sigma_{ess}(M\setminus K_{i})$.

 By hypothesis we have a complete  $m$-submanifold $\varphi\colon M\hookrightarrow B_{N}(r)$ properly immersed into a ball $B_{N}(r)=B_{N}(p,r)$ with center at $p$ and radius $r$ in a Riemannian $n$-manifold $N$ with radial sectional
curvatures $K_N^{\mathrm{rad}}$ along the radial geodesics issuing from
$p$  bounded
as $a=\inf K_{N}^{\mathrm{rad}}\leq K_N^{\mathrm{rad}}\leq b=\sup K_{N}^{\mathrm{rad}}$ in $B_N(r)$, where  $r<\min\{{\rm inj}_N(p), \pi/2\sqrt{b}\}$. Here we replace $\pi/2\sqrt{b}$ by $+\infty$ if $b\leq 0$.

\vspace{2mm}
Define a function  $v\colon B_{N}(p,r)\to \mathbb{R}$  by $v(y)=\phi_{a}(\rho (y))$, where  $\phi_{a}:[0,r]\to \mathbb{R}$ given by
\begin{equation}\phi_{a} (t)=\left\{ \begin{array}{lll}
\cos(\sqrt{a}\,t)-\cos(\sqrt{a}\, r) & \mathrm{if} & a>0, \;\; t<\pi/2\sqrt{a},\\
r^{2}-t^{2} & \mathrm{if} &  a=0,\\
\cosh(\sqrt{-a}\, r)-\cosh(\sqrt{-a}\, t) & \mathrm{if} & a<0.
\end{array}\right.\label{eqphi}
\end{equation}
 Observe that $\phi (t)>0$ in $[0, r)$, $\phi_{a}(r)=0$, $\phi_{a}'(t)<0$ and $\phi_{a}''(t)-C_{a}(t)\phi_{a}'(t)=0$ in $[0,r]$.  This function $\phi_{a}$ we learned from Markvorsen \cite{markvorsen2}. Let $f\colon M\to \mathbb{R}$ defined by $f=v\circ \varphi$ and consider an exhaustion sequence of $M$ by compact sets  $K_{i}=\varphi^{-1}\left(\overline{B_{N}(p,r_{i})}\right)$, where  $r_{i}< r$, $r_{i}\to r$.
  By Barta's Theorem  we have that $\lambda^{\ast}(M\setminus K_{i})\geq \inf_{M\setminus K_{i}}( \displaystyle\frac{-\triangle f}{f}).$

\vspace{2mm}
\noindent Now by (\ref{eqBF3}) we have
\begin{eqnarray}\triangle f (x)&=& \sum_{i=1}^{m}\hess_{N}v(\varphi (x))(e_{i},e_{i})+\langle \grad\,v\,,\, \sum_{i=1}^{m}\alpha
(e_{i},e_{i})\rangle.\nonumber \\
&=&\sum_{i=1}^{m}\hess_{N}v(\varphi (y))(e_{i},e_{i})+\langle \grad\,v\,,\, H\rangle.\nonumber
\end{eqnarray}
The metric of $N$ inside the normal geodesic ball $B_{N}(p,r)$ can be written in polar coordinates as $ds^{2}=dt^{2}+\vert A(t,\xi)\vert^{2}d\xi^{2}$,  ${\mathcal A}(t,\xi)$ satisfies the Jacobi
equation ${\mathcal A}''+{\mathcal R}{\mathcal A}=0$ with initial conditions
${\mathcal A}(0,\xi)=0$, ${\mathcal A}'(0,\xi)=I$.  We have at the point $\varphi (x)$ an orthonormal basis $\{\partial/\partial t, \partial/\partial \xi_{1}, \ldots \partial/\partial \xi_{n-1}\}$ for $T_{\varphi(x)}N$. We may choose an orthonormal basis for $T_{x}(M\setminus K_{i})$ as  $e_{1}=\langle e_{1}, \partial/\partial t\rangle \cdot \partial/\partial t + e_{1}^{\perp}$, where $ e_{1}^{\perp}\perp \partial/\partial t$ and $\{e_{2}, \ldots, e_{m}\}\subset \{\partial/\partial \xi_{1}, \ldots \partial/\partial \xi_{n-1}\}$. Computing $\hess_{N}v(\varphi (x))(e_{i},e_{i})$ we have

 \begin{eqnarray} \hess_{N}v(\varphi (x))(e_{1},e_{1})&=&\left[\phi_{a}''(t)-\phi_{a}'(t)\cdot \hess_{N} \rho (e_{1}^{\perp}/\vert e_{1}^{\perp}\vert ,e_{1}^{\perp}/\vert e_{1}^{\perp}\vert)\right]\langle e_{1}, \grad \rho\rangle^{2}\nonumber \\ &&   +\,\, \phi_{a}'(t)\cdot \hess_{N} \rho (e_{1}^{\perp}/\vert e_{1}^{\perp}\vert) ,e_{1}^{\perp}/\vert e_{1}^{\perp}\vert))\nonumber \\
 &=& \phi_{a}'(t)\cdot \left[C_{a}(t)-\hess \rho (e_{1}^{\perp}/\vert e_{1}^{\perp}\vert ,e_{1}^{\perp}/\vert e_{1}^{\perp}\vert)\right]\langle e_{1}, \grad \rho\rangle^{2} \\
  && +\,\, \phi_{a}'(t)\cdot \hess_{N} \rho (e_{1}^{\perp} ,e_{1}^{\perp})\nonumber
 \end{eqnarray}
and for $i\geq 2$
\begin{equation}\hess_{N}v(\varphi (x))(e_{i},e_{i})= \displaystyle \phi_{a}'(t)\cdot \hess_{N} \rho (e_{i}, e_{i}) \end{equation}
 where $t=\rho( \varphi (x))$.  Now,
 \begin{eqnarray}-\triangle f &= &- \phi_{a}'(t)\cdot \left[C_{a}(t)-\hess \rho (e_{1}^{\perp}/\vert e_{1}^{\perp}\vert) ,e_{1}^{\perp}/\vert e_{1}^{\perp}\vert))\right]\langle e_{1}, \grad \rho\rangle^{2} \nonumber \\
  && -\,\, \phi_{a}'(t)\cdot\left[ \hess_{N} \rho (e_{1}^{\perp}/\vert e_{1}^{\perp}\vert) ,e_{1}^{\perp}/\vert e_{1}^{\perp}\vert))+\sum_{i=2}^{m} \hess_{N} \rho (e_{i}, e_{i})\right]\nonumber\\
  && -\,\,\phi_{a}'(t)\langle \grad \rho, H\rangle \nonumber \\
  &\geq & -\phi_{a}'(t)\cdot \left[ m\cdot C_{b}(t)-\sup \vert H\vert \right]\nonumber
 \end{eqnarray} We used that  $C_{a}(t)\geq \hess \rho (e_{1}^{\perp}/\vert e_{1}^{\perp}\vert) ,e_{1}^{\perp}/\vert e_{1}^{\perp}\vert))\geq C_{b}(t)$, $\hess \rho (e_{i}, e_{i})\geq C_{b}(t)$ by the Hessian Comparison Theorem and  that  $-\phi_{a}'(t)>0$.

 \vspace{3mm}
 Hence \begin{eqnarray}\lambda^{\ast}(M\setminus K_{i})\geq \displaystyle\inf_{M\setminus K_{i}} ( -\displaystyle\frac{\triangle f}{f})& \geq & \inf_{M\setminus K_{i}}-\displaystyle\frac{\phi_{a}'(t)}{\phi_{a}(t)}\left[ m\cdot C_{b}(t)-\sup \vert H\vert \right]\nonumber \\ &\geq & -\displaystyle\frac{\phi_{a}'(r_{i})}{\phi_{a}(r_{i})}\left[ m\cdot C_{b}(r)-\sup \vert H\vert \right]
 \end{eqnarray} Thus $\lambda^{\ast}(M\setminus K_{i})\to \infty$ as $r_{i}\to r$.  This proves Theorem (\ref{thmMain}).
\section{cylindrically bounded submanifolds.}
Let $\varphi \colon M^{m}\hookrightarrow B_{N}(r)\times \mathbb{R}^{\ell}\subset N^{n-\ell}\times \mathbb{R}^{\ell}$, $m\geq \ell+1$,  be an isometric immersion of a complete  Riemannian $m$-manifold $M^{m}$ into the $B_{N}(r)\times \mathbb{R}^{\ell}$, where $B_{N}(r)$ is a geodesic ball in a Riemannian ($n-\ell$)-manifold $N^{n-\ell}$, centered at a point $p$ with radius $r$.  Let  $b=\sup K_N^{\mathrm{rad}}$ where $K_N^{\mathrm{rad}} $ are  the radial sectional
curvatures along the  geodesics issuing from
$p$. Assume that  $r<\min\{{\rm inj}_N(p), \pi/2\sqrt{b}\}$, where $\pi/2\sqrt{b}=+\infty$ if $b\leq 0$.
\begin{theorem}\label{thmMain2}Suppose that $\varphi \colon M^{m}\hookrightarrow B_{N}(r)\times \mathbb{R}^{\ell}$ as above satisfies the following.\begin{itemize}\item[1.] For every $s<r$, the set $\varphi^{-1}(B_{N}(s)\times \mathbb{R}^{\ell})$ is compact in $M$.
\item[2.] $\sup_{M}\vert H\vert<(m-\ell)C_{b}(r)$ where $\vert H\vert(x) $ is norm of the mean curvature vector of $\varphi (M) $ at $\varphi (x)$.
 \end{itemize}Then $M$ has discrete spectrum.
\end{theorem}
Observe that the condition 1. is a stronger property than  being a proper immersion except when $\ell=0$.
 \vspace{2mm}

\noindent {\em Proof:} As before, let
 $r_{i}\to r$ be a sequence of positive real numbers  $r_{i}<r$ and the compacts sets   $K_{i}=\varphi^{-1}(B_{N}(r_{i})\times \mathbb{R}^{\ell})$. We need only to show that $\lambda^{\ast}(M\setminus K_{i})\to \infty $ as $r_{i}\to r$.  Define $v$ on $B_{N}(r)\times \mathbb{R}^{\ell}$ by $v(x,y)=\phi_{a}(\rho(x))$, where $\rho(x)={\rm dist_{N}}(p,x)$, $\phi_{a}$ given in (\ref{eqphi}) and $a=\inf K_{N}^{\mathrm{rad}}$. Let $f=v\circ \varphi:M\to \mathbb{R} $.
We have by (\ref{eqBF3})
\begin{eqnarray} \triangle f (x)& = &  \sum_{i=1}^{m}\hess_{N\times \mathbb{R}^{\ell}}v(\varphi (x))(e_{i},e_{i})+  \langle \grad v,\, H\rangle.\nonumber \\
&=& \sum_{i=1}^{m}\hess_{N}(\phi_{a}\circ \rho) (\varphi (x))(e_{i},e_{i})+  \langle \grad (\phi_{a}\circ \rho) ,\, H\rangle.\end{eqnarray}
 At  $\varphi (x)=(y_{1},y_{2})$, consider the orthonormal  basis  $$\{\stackrel{Polar\, basis}{\overbrace{\partial/\partial t, \partial/\partial \xi_{1}, \ldots \partial/\partial \xi_{n-\ell-1}}},\stackrel{Cartesian\, basis}{\overbrace{ \partial/\partial s_{1}, \ldots \partial/\partial s_{\ell}}}\}$$ for $T_{(y_{1},y_{2})}N^{n-\ell}\times \mathbb{R}^{\ell}=T_{y_{1}}N^{n-\ell}\oplus T_{y_{2}}\mathbb{R}^{\ell}$. Choose an orthonormal basis $\{e_{1}, e_{2}, \ldots, e_{m}\}$ as follows

\begin{equation*} e_{i} = a_{i}\frac{\partial}{\partial t}  + \displaystyle \sum_{j=1}^{n-\ell-1} b_{ij} \frac{\partial }{\partial \xi_{j}} +
 \displaystyle \sum_{j=1}^{\ell} c_{ij} \frac{\partial}{\partial s_{j}}. \end{equation*}

Using that $\phi_{a}'(t)<0$, $\phi_{a}''=C_{a}(t)\phi_{a}'(t)$ and $\hess \rho(y_{1}) (\partial/ \partial \xi_{j}, \partial /\partial \xi_{j})\geq C_{b}(t)$ for all $j=1,\ldots, n-\ell-1$, we have that \begin{eqnarray*}\hess \phi_{a}\circ\rho (y_{1})(e_{i},e_{i})&=& \phi_{a}''(t)a_{i}^{2} + \phi_{a}'(t) \sum_{j=2}^{n-\ell-1} b_{ij}^{2}\,\hess \rho(y_{1}) (\partial/ \partial \xi_{j}, \partial /\partial \xi_{j})\nonumber\\
&\leq & \phi_{a}''(t)a_{i}^{2} + \phi_{a}'(t) \sum_{j=2}^{n-\ell-1} b_{ij}^{2}\,C_{b}(t)\nonumber \\
&=& C_{a}(t) \phi_{a}'(t)a_{i}^{2} + \phi_{a}'(t)(1-a_{i}^{2}-\sum_{k=1}^{\ell} c_{ik}^{2})C_{b}(t)\\
&=& \phi_{a}'(t)a_{i}^{2}(C_{a}(t)-C_{b}(t))+ \phi_{a}'(t)(1-\sum_{k=1}^{\ell} c_{ik}^{2})C_{b}(t)\\&\leq &\phi_{a}'(t)(1-\sum_{k=1}^{\ell} c_{ik}^{2})C_{b}(t)
\end{eqnarray*} since $C_{a}(t)\geq C_{b}(t)$ and where $t=\rho (y_{1})$. Therefore
\begin{eqnarray*}-\sum_{i=1}^{n}\hess \phi_{a}\circ\rho (y_{1})(e_{i},e_{i})&\geq & -\phi_{a}'(t)(m-\sum_{i=1}^{m}\sum_{k=1}^{\ell} c_{ik}^{2})C_{b}(t)\\
&\geq & -\phi_{a}'(t)(m-\ell)C_{b}(t)
\end{eqnarray*}
From this we have that
\begin{eqnarray}-\displaystyle \frac{\triangle f }{f}(x)&\geq & -\frac{\phi_{a}'(t)}{\phi_{a}(t)}\left[ (m-\ell)C_{b}(t) - \sup_{M} \vert H\vert\right]
\end{eqnarray}
 so that \begin{eqnarray*}\displaystyle \inf_{M\setminus K_{i}}(-\displaystyle \frac{\triangle f }{f}) &\geq& -\frac{\phi_{a}'(r_{i})}{\phi_{a}(r_{i})}\left[ (m-\ell)C_{b}(r) - \sup_{M} \vert H\vert\right].
\end{eqnarray*}
Therefore $ \displaystyle \inf_{M\setminus K_{i}}(-\displaystyle \frac{\triangle f }{f})\to +\infty$ as $r_{i}\to r$ proving Theorem (\ref{thmMain2}).

\bibliographystyle{amsplain}

\end{document}